\documentclass{article}
\usepackage{a4wide}
\usepackage{amssymb}
\usepackage{amsfonts}
\usepackage{amsxtra}
\usepackage{latexsym}
\usepackage{graphicx}
\usepackage[colorlinks=true]{hyperref}
\usepackage[capitalise]{cleveref}
\usepackage{authblk}
\usepackage{amsmath}

\begin{document}

\newtheorem{theo}{Theorem}[section]
\newtheorem{theo*}{Theorem}
\newtheorem{lemma}[theo]{Lemma}
\newtheorem{lemma*}{Lemma}
\newtheorem{definition}[theo]{Definition}
\newtheorem{notation}[theo]{Notation}
\newtheorem{obs}[theo]{Observation}
\newtheorem{remark}[theo]{Remark}
\newtheorem{cor}[theo]{Corollary}
\newtheorem{cor*}{Corollary}
\newtheorem{prop}[theo]{Proposition}
\newtheorem{conj}[theo]{Conjecture}
\newtheorem{example}[theo]{Example}

\newcommand{\N}{{\mathbb{N}}}
\newcommand{\Z}{{\mathbb{Z}}}   
\newcommand{\Q}{{\mathbb{Q}}}
\newcommand{\R}{{\mathbb{R}}}
\newcommand{\C}{{\mathbb{C}}}
\newcommand{\M}{{\mathbb{M}}}
\newcommand{\F}{{\mathbb{F}}}
\newcommand{\V}{{\mathbb{V}}}
\newcommand{\LL}{{\mathbb{L}}}
\newcommand{\A}{{\mathbb{A}}}
\newcommand{\HH}{{\mathbb{H}}}

\newcommand{\g}[1]{{\mathfrak {#1}}}
\newcommand{\cc}[1]{{\cal {#1}}}
\newcommand{\ul}{\underline }
\newcommand{\ol}{\overline }

\newcommand{\qed}{\hspace*{\fill}$\Box$}
\newcommand{\Qed}{\hspace*{\fill}$\Box \Box$}


\title{On the first-order theory of the remainder}

\author{Mihai Prunescu\thanks{Research Center for Logic, Optimization and Security (LOS), Faculty of Mathematics and Computer Science, University of Bucharest, Academiei 14, 010014 Bucharest, Romania; and  Simion Stoilow Institute of Mathematics of the Romanian Academy, Research unit 5, P. O. Box 1-764, RO-014700 Bucharest, Romania.
mihai.prunescu@imar.ro, mihai.prunescu@gmail.com}
}

\date{}
\maketitle

\begin{abstract}
\parindent 0 cm  
It is proved that the first-order theory of the structure $(\mathbb N, \bmod)$ is undecidable. Here $\bmod $ denotes the operation of computing the remainder for any division between positive integers; i.e, $x \bmod y$ is the remainder obtained by the division $x : y$.  
\vspace{1mm}


\vspace{1mm}

\thanks{A.M.S.-Classification: 03D35 }
\end{abstract} 


Denote with $\mathbb N$ the set of positive integers together with $0$ and with $\mathbb N_+$ the set of positive integers.
By definition, for $x \in \mathbb N$ and $y \in \mathbb N_+$, the number $x \bmod y$ is the unique natural number such that $0 \leq x \bmod y < y$ and there is a $q \in \mathbb N$ with $x = qy+ (x \bmod y)$. 
So, $5 \bmod 2 = 1$, $5 \bmod 3 = 2$, $6 \bmod 3 = 0$, etc. 

A classical example of decidable fragment of the arithmetic is Presburger's Arithmetic, i.e. the first-order theory of the structure $(\mathbb N, +)$, see Presburger \cite{Presburger}. To prove the decidability, one expands the structure by defining several relations, including the total order $<$ and the infinite family of congruence relations $(a \equiv b \mod n)$ with $n \geq 2$ and one shows that the expanded structure allows an effective method of quantifier elimination according to the new language. 

If instead of the infinite family of congruence relations, we consider the remainder operation $\bmod$, the situation changes.  We show here that this structure has an undecidable theory. 

As we need a totally defined operation,  we choose some artificial definition for the values of $n \bmod 0$, say:
$$\forall \, n \in \mathbb N \,\,\,\, n \bmod 0 := n,$$
including $0 \bmod 0 = 0$. 

We will make a connection between the structure $(\mathbb N, \bmod)$ and two classical results of Julia Robinson, see \cite{Robinson}. We express this connection in a slightly more general form as follows:

\begin{lemma*}\label{lemma}
We recall that $y \mid x$ is the relation of divisibility and $x < y$ the strict order relation.
Let $ \# : \mathbb N^2 \rightarrow \mathbb N$ be an operation with the following properties. For all $x, y \in \mathbb N$ one has: 
\begin{enumerate} 
\item $x = 0$ if and only if $\exists \,y\,\, y \# y = x$. Observe that this is equivalent with  $$(\mathbb N,  \# ) \models \forall \,y\,\,\,\,y \# y = 0.$$ 
\item $x < y$ if and only if $( x = 0 \,\wedge \, y \neq 0) \vee ( x \neq 0 \wedge y \neq 0 \wedge x \# y = x)$. 
\item $y \mid x$ if and only if $x = 0 \vee ( x \neq 0 \wedge x \# y = 0)$. 
\end{enumerate} 
Then the first-order theory of the structure $(\mathbb N, \# )$ is undecidable. 
\end{lemma*}

{\bf Proof}: According to Theorem 1.1 of Robinson \cite{Robinson}, in the set $\mathbb N_+$ of positive natural numbers, the addition is first-order definable in terms of multiplication and the successor operation, as well as in terms of multiplication and the relation of less-than. Also, according with Theorem 1.2 of the same paper, the multiplication of positive integers is arithmetically definable in terms of the successor operation and the relation of divisibility. But given the total order $<$, the successor operation $S(x) := x+1$ is definable by:
$$ x = S(y) \textrm{ if } x > y \,\wedge \, \forall \,z\,( z > x \rightarrow ( z = y \,\vee\, z > y )).$$
It follows that in the structure of positive natural numbers $(\mathbb N_+, <, \mid)$, both the addition and the multiplication are first-order definable. According to Tarski, Mostowski and Robinson  \cite{TMR}, the first-order theory of the structure $(\mathbb N_+, <, \mid)$ is undecidable. 

Now, we use the properties of the operation $ \#$ to interpret the structure $(\mathbb N_+, <, \mid)$ in the structure $(\mathbb N, \# )$. The universe of the structure is simply defined by:
$$x \in \mathbb N_+ \textrm{ if } x \neq 0,$$
where we apply the fact that $x = 0$ is a definable relation in $(\mathbb N,  \# )$. The order and the divisibility for elements $a, b \in \mathbb N_+$ are defined by:
$$a < b \textrm{ if } a \# b = a,$$
$$\,\,\, b \mid a \textrm{ if } a \# b = 0.$$
The second definition uses again the fact that $0$ is definable in $(\mathbb N,  \# )$. We also observe that the relations of order and divisibility defined above on $\mathbb N_+$ are the intersections of the corresponding subsets of  $\mathbb N \times \mathbb N$ called $<$, and respectively $\mid$,  with the set $\mathbb N_+ \times \mathbb N_+$ - so they are natural restrictions of these relations. 

By the usual interpretation trick, this leads to an algorithmic embedding of the complete first-order theory of the structure $(\mathbb N_+, <, \mid)$ in the complete first-order theory of the structure $(\mathbb N, \# )$, so the latter is undecidable. \qed 

Now we can state the main result:

\begin{cor*}\label{main}
The first-order theory of the remainder operation $(\mathbb N, \bmod)$ is undecidable. 
\end{cor*} 

{\bf Proof}: It is immediate to check that the operation $\# = \bmod : \mathbb N^2 \rightarrow \mathbb N$  satisfies the conditions of  Lemma \ref{lemma}. \qed

\end{document}